\documentclass{amsproc}
\usepackage{graphicx}
\usepackage{amscd}
\usepackage{amsmath}
\usepackage{amsfonts}
\usepackage{amssymb}
\newtheorem{theorem}{Theorem}
\theoremstyle{plain}

\newtheorem{corollary}{Corollary}
\newtheorem{definition}{Definition}

\numberwithin{equation}{section}

\begin{document}
\title[Peetre's Interpolation Spaces]{On Properties of Symmetric Banach Function spaces and Peetre's Interpolation Spaces}
\author{Eugene Tokarev}
\address{B.E. Ukrecolan, 33-81 Iskrinskaya str., 61005, Kharkiv-5, Ukraine}
\email{evtokarev@univer.kharkov.ua}
\subjclass{Primary 46E30, 46B70; Secondary 46B10, 46B20}
\keywords{Symmetric Banach spaces of measurable functions, Interpolation spaces, Partial
order, Lattices}
\dedicatory{Dedicated to the memory of S. Banach.}
\begin{abstract}A set $\mathcal{S}$ of all symmetric Banach function spaces defined on [0,1]
is equipped with a partial order by the relation $\subset^{c}$ of continuous
inclusion. Properties of symmetric spaces, which does not depend of their
position in the ordered structure are studied. With the help of the J.
Peetre's interpolation scheme it is shown that for any pair $E$, $F$ of
symmetric spaces such that $F$ is absolutely continuously included in $E$
there exists an intermediate (so called, Peetre's) space $K(E,F;W)$, where $W$
is a space with an unconditional basis, that is reflexive or, respectively,
weakly sequentially complete if $W$ also is reflexive or, respectively, weakly
sequentially complete.
\end{abstract}
\maketitle

\section{ Introduction}

A Banach space $E$ of (classes of) Lebesgue measurable real functions acting
on the interval $[0,1]$ (equipped with the Lebesgue measure
$\operatorname*{mes}(e)$ for a measurable subset $e\subset\lbrack0,1]$) is
called \textit{symmetric }if for any pair of (measurable) functions
$x=x\left(  t\right)  $ and $y=y\left(  t\right)  $ following conditions are satisfied:

1\textit{. If }$x\in E$\textit{ and }$\left|  y(t)\right|  \leq\left|
x(t)\right|  $\textit{ (a.e) then }$y\in E$\textit{ and }$\left\|  y\right\|
_{E}\leq\left\|  x\right\|  _{E}$\textit{.}

\textit{2. If }$x\in E$\textit{ and }$\left|  y(t)\right|  $\textit{ is
equimeasurable with }$\left|  x(t)\right|  $\textit{ then }$y\in E$\textit{
and }$\left\|  y\right\|  _{E}=\left\|  x\right\|  _{E}$.

Let $\mathcal{S}$ be a class of all \textit{symmetric Banach function spaces}
(which in the future will be referred to as \textit{symmetric spaces}). The
class $\mathcal{S}$ may be partially ordered by a relation $E\subset^{c}F$,
which means that $E$, as a set of functions, is contained in $F$ and a norm
$\left\|  x\right\|  _{E}$ of any function $x\in E$ is estimated by $\left\|
x\right\|  _{E}\leq c(E,F)\left\|  x\right\|  _{F}$. To be more precise,
define an operator $i\left(  E,F\right)  :E\rightarrow F$, which asserts to
any $x\in E$ the same function $x\in F$. This operator is called \textit{the
inclusion operator}. Its norm (which is equal to the infimum of all possible
constants $c(E,F)$) is called \textit{the inclusion constant}.

A lot of properties of symmetric spaces depend on their position in a
partially ordered set $\left\langle \mathcal{S},\subset^{c}\right\rangle $.
E.g., extreme spaces, i.e., the minimal and the maximal elements of
$\left\langle \mathcal{S},\subset^{c}\right\rangle $, are exactly classical
spaces $L_{\infty}[0,1]$ and $L_{1}[0,1]$ respectively (cf. [1]). Let $G$
denotes an Orlicz space, which norming function is given by the equality
$M(u)=(\exp(u^{2})-1)/(e-1)$, and $G_{0}$ denotes a closure of $L_{\infty
}[0,1]$ in $G$. Any symmetric space $E$ such that $G_{0}\subset^{c}%
E\subset^{c}\left(  G_{0}\right)  ^{\ast}$ has a basis of functions, which are
uniformly bounded in the $L_{\infty}$-norm (cf. [2]). A lot of other similar
examples may be found in [3], [4], [5].

Nevertheless there are properties of symmetric spaces $E$, which do not depend
on position of $E$ in the structure $\left\langle \mathcal{S},\subset
^{c}\right\rangle $. Some of such properties are studied in the article. In
the study it will be used an interpolation construction due to J. Peetre [6],
which asserts to a pair $\left(  E,F\right)  $ of symmetric spaces with
$F\subset^{c}E$ a new Banach space $K(E,F;W)$ (its construction uses a given
Banach space $W$ with an unconditional basis) that has certain interpolation
properties with respect to spaces $E$ and $F.$

This construction (and its generalizations) plays an important role in the
functional analysis: in the interpolation theory (cf. [7]), in factorizations
of weakly compact operators (cf. [8]), in constructing of peculiar examples of
Banach spaces ([9], [10]) and many others.

The construction of the space $K(E,F;W)$ may one led to believe that its
properties hardly depend on properties of $E$ and of $F$.

However, if an order interval $[F,E]_{\subset^{c}}$ is \textit{sufficiently
large} (exactly: if $F$ \textit{is absolutely continuously included in }$E$;
shortly: $F\Subset E$; this relation will be explained later) then some
properties of $K(E,F;W)$ depend only on properties of $W$. These properties
will be called $\omega$-properties and will be defined below. In particular,
such is a property of a symmetric space to be either reflexive or weakly
sequentially complete.

\section{Definitions and notations}

Let $E$ be a symmetric space. Its \textit{fundamental function} $\varphi
_{E}\left(  \tau\right)  $ is given by
\[
\varphi_{E}\left(  \tau\right)  =\left\|  \chi_{\lbrack0,\tau]}\left(
t\right)  \right\|  _{E},
\]
where $\chi_{e}\left(  t\right)  $ is an \textit{indicator function} of a
measurable subset $e\subset\lbrack0,1]$:
\[
t\in e\Rightarrow\chi_{e}\left(  t\right)  =1;t\notin e\Rightarrow\chi
_{e}\left(  t\right)  =0
\]
In the future it will be assumed that $\varphi_{E}\left(  0\right)  =0$;
$\varphi_{E}\left(  1\right)  =1$.

$E$ has an \textit{absolutely continuous norm} provided
\[
\lim\nolimits_{\operatorname*{mes}(e)\rightarrow0}(\left\|  x(t)\chi
_{e}\left(  t\right)  \right\|  _{E}=0
\]
for any $x(t)\in E$. It is known (see e.g. [11]) that $E$ has such a norm if
and only if it is separable.

A sequence $\{f_{n}:n<\infty\}$ may be also referred to as $(f_{n})_{n<\infty
}$ or, simply, as\ $(f_{n})$. A sequence $(f_{n})\subset E$ is said to
be\textit{\ disjoint} if $\ f_{m}(t)f_{n}(t)=0$ (a.e.) for $n\neq m$.

Let $E$ be a symmetric space, $x(t)\in E$, $\left\|  x\right\|  _{E}\neq0$,
$B\subset E$. Define following characteristics:
\[
\eta_{E}(x,\tau)=\sup\{\left\|  x(t)\chi_{e}(t)\right\|  _{E}/\left\|
x\right\|  _{E}:\operatorname*{mes}(e)=\tau\};
\]%
\[
\eta_{E}(B,\tau)=\sup\{\eta_{E}(x,\tau):x\in B;\text{ \ }\left\|  x\right\|
_{E}\neq0\};
\]%
\[
\eta_{E}(B)=\underset{\tau\rightarrow0+}{\lim}\eta_{E}(B,\tau).
\]
In a case when $E=L_{1}[0,1]$ the lower index $_{E}$ will be omitted. Thus,
$\eta(x,\tau)$ denotes $\eta_{L_{1}}(x,\tau)$ and so on.

Let $E$, $F$ be a pair of symmetric spaces; $F\subset^{c}E$; $i:F\rightarrow
E$ be the \textit{inclusion operator}. It will be said that $F$ is
\textit{absolutely continuously included in} $E$, in symbols: $F\Subset E$, if
$\eta_{E}(iF)=0$.

Let $a$, $b$ be positive real numbers; $E$, $F$ be symmetric spaces; $f$ be a
measurable function. Put
\[
k(f;a,b)=\inf\{a\left\|  u\right\|  _{E}+b\left\|  v\right\|  _{F}:u+v=f\}.
\]
Clearly, $k(f;a,b)$ may be regarded as a norm on a set $E+F$ of all functions
$f$ of kind $f=u+v$, for which $k(f;a,b)$ is finite.

\begin{definition}
Let $E$, $F$ be symmetric spaces; $F\Subset E$; $W$ be a Banach space with an
unconditional basis $(w_{i})$. Let $\left(  a_{i}\right)  $ and $\left(
b_{i}\right)  $ be sequences of positive real numbers such that $\left(
a_{i}\right)  $ is non decreasing; $\left(  b_{i}\right)  $ is non increasing;
$\lim_{i\rightarrow\infty}a_{i}=\infty$; $\sum\nolimits_{i<\infty}b_{i}%
<\infty$.

A Peetre's space $K(E,F;W;(a_{i}),(b_{i}))$ is a Banach space of measurable
functions $f(t)$, defined on $[0,1]$, such that a series $\sum
\nolimits_{i<\infty}k(f;a_{i},b_{i})w_{i}$ convergents in $W$. A norm on the
Peetre's space (that will be below referred to as $K(E,F;W)$ or, simply, as
$K$) is given by
\[
\left\|  f\right\|  _{K}=\left\|  \sum\nolimits_{i<\infty}k(f;a_{i}%
,b_{i})w_{i}\right\|  _{W}/\left\|  \chi_{\lbrack0,1]}(t)\right\|  _{K}.
\]
\end{definition}

Obviously, $K(E,F;W)$ is a symmetric space with an absolutely continuous norm
(and, hence, is separable). As it was mentioned before, it will be supposed
that $\varphi_{K}(1)=\left\|  \chi_{\lbrack0,1]}(t)\right\|  _{K}=1$.

Let $(x_{n})\subset X$ be a sequence of elements of a Banach space $X$. A
sequence $(y_{k})\subset X$ is said to be a \textit{block sequence with
respect to }$(x_{n})$ if there exist a sequence $(a_{n})$ of reals and a
sequence $1=n_{0}<n_{1}<...<n_{k}<...$ of naturals such that for all
$k<\infty$%
\[
y_{k}=\sum\nolimits_{m=n_{k}-1}^{n_{k+1}}a_{m}x_{m}%
\]
If $(x_{n})$ is a basic sequence then a block sequence $(y_{k})$ with respect
to $(x_{n})$ will be called a \textit{block-basis sequence.}

\section{$\omega$-properties of Peetre's spaces.}

A property (P) of a Banach space $X$ is said to be an $\omega$%
\textit{-property} if a condition ''\textit{every subspace of }$X$%
\textit{\ has a subspace that has the property (P)}'' implies that $X$ also
has the (P). The reflexivity (R) and the weak sequential completeness (WSC)
are not $\omega$-properties of general Banach spaces (take a look, e.g., at
the classical James' quasireflexive space $J$). Nevertheless on the class
$\mathcal{S}$ of symmetric spaces both of them, (R) and (WSC), are $\omega
$-properties. Indeed, if $E\in\mathcal{S}$ is not reflexive then it contains a
subspace, that is isomorphic either to $l_{1}$ or to $c_{0}$; if $E$ is not
weakly sequentially complete then $E$ contains a subspace which is isomorphic
to $c_{0}$. Clearly, every subspace of both $l_{1}$ and $c_{0}$ is
non-reflexive; every subspace of $c_{0}$ fails to have the property (WSC).

\begin{theorem}
Let $E$, $F$ be symmetric spaces; $F\Subset E$; $(a_{i})$; $(b_{i})$ be as in
definition 1; $W$ be a space with an unconditional basis $(w_{i})$. Then every
disjoint sequence of functions $(f_{n})\subset K(E,F;W)$ contains a
subsequence $(f_{n_{m}})$, which is equivalent to a some block-basic sequence
with respect to $(w_{n}).$
\end{theorem}

\begin{proof}
Let $x(t)$ be a measurable function. Put $s(\tau)=\eta_{E}(F,\tau)$. Since
$F\Subset E$, it is clear that $\underset{\tau\rightarrow0+}{\lim}s(\tau)=0$.
Let
\[
A_{n}=\sum\nolimits_{k=1}^{n}a_{k}.
\]
Since $a_{k}\uparrow\infty$, for every $\varepsilon_{0}>0$ and every $n_{0}%
\in\mathbb{N}$ there exists a number $j=j(\varepsilon_{0},n_{0})$ such that
$2A_{n}a_{j}^{-1}<\varepsilon_{0}$. Fix $n_{0}$ and $\varepsilon_{0}$ and
consider a function $f(t)\in K$; $\left\|  f\right\|  _{K}=1$\ with a support
$A_{f}=\operatorname*{supp}(f):=\{t\in\lbrack0,1]:f(t)\neq0\}.$ It will be
assumed that the support of $f$ \ is sufficiently small, namely:
$s(\operatorname*{mes}(A_{f}))\leq a_{j}^{-1}b_{j}$, where $j=j(\varepsilon
_{0},n_{0})$ have been chosen before. It will be shown that a norm of the
function $f$ can be calculated up to $\varepsilon_{0}$ only by using vectors
$w_{n}$ with $n>n_{0}$, and that $n_{0}$ tends to infinity provided
$\operatorname*{mes}(\operatorname*{supp}(f))$ tends to zero (this procedure
will be called '\textit{a} \textit{cutting of the} \textit{head}'). Let
$\left\|  f\right\|  _{K}=1$. Then
\[
\inf\left\{  a_{j}\left\|  u\right\|  _{E}+b_{j}\left\|  v\right\|
_{F}:u+v=f\right\}  \leq1.
\]
Thus, $\left\|  u\right\|  _{E}\leq a_{j}^{-1}$; $\left\|  v\right\|  _{F}\leq
b_{j}^{-1}.$

Without loss of generality it may be assumed that
\[
\operatorname*{supp}(u)\subset\operatorname*{supp}(f);\operatorname*{supp}%
(v)\subset\operatorname*{supp}(f).
\]
Indeed, in a contrary case from the equality $u\chi_{A_{f}}+v\chi_{A_{f}%
}=f\chi_{A_{f}}=f$ and because of the symmetry, it follows:
\[
\left\|  u\chi_{A_{f}}\right\|  _{E}\leq\left\|  u\right\|  _{E};\text{
\ \ \ }\left\|  v\chi_{A_{f}}\right\|  _{F}\leq\left\|  v\right\|  _{F}.
\]
Clearly, $\left\|  f-v\right\|  _{E}=\left\|  u\right\|  _{E}\leq a_{j}^{-1}$.
From the triangle inequality it follows that
\begin{align*}
\left\|  f\right\|  _{E} &  \leq\left\|  u\right\|  _{E}+\left\|  v\right\|
_{E}\leq a_{j}^{-1}+\left\|  v\right\|  _{E}\\
&  \leq a_{j}^{-1}+b_{j}^{-1}\eta_{K}(v,\operatorname*{mes}(A_{f}))\leq
a_{j}^{-1}+b_{j}^{-1}a_{j}^{-1}b_{j}\leq2a_{j}^{-1}.
\end{align*}
Since $F\Subset E$,
\[
\max(a_{n},b_{n})k(f;1,1)\geq k(f;a_{n},b_{n}).
\]
Hence,
\[
a_{m}^{-1}k(f;a_{m},b_{m})\leq\left\|  f\right\|  _{E}\leq2a_{j}^{-1}.
\]
Consequently,
\[
k(f;a_{m},b_{m})\leq2a_{j}^{-1}a_{m}.
\]
With no loss of generality it may be assumed that the basis $\left(
w_{n}\right)  $ is 1-unconditional, i.e. that $\left\|  \sum a_{i}%
w_{i}\right\|  \leq\sum\left|  a_{i}\right|  $ for any choice of scalars
$\left(  a_{i}\right)  $. So,
\begin{align*}
\left\|  \sum\nolimits_{m=1}^{n}k(f;a_{m},b_{m})w_{m}\right\|  _{W} &
\leq\left\|  \sum\nolimits_{m=1}^{n}2a_{j}^{-1}a_{m}w_{m}\right\|  _{W}\\
&  \leq2(\sum\nolimits_{m=1}^{n}a_{m})a_{j}^{-1}=2A_{n}a_{j}^{-1}%
\leq\varepsilon_{0}.
\end{align*}
The series $\sum_{m<\infty}k(f;a_{m},b_{m})w_{m}$ convergents in $W$. Hence,
there exists $N=N(\varepsilon_{0})$ such that
\[
\left\|  \sum\nolimits_{N<m<\infty}k(f;a_{m},b_{m})w_{m}\right\|  _{W}%
\leq\varepsilon_{0}.
\]
(\textit{''a cutting of the tail''}). Therefore,
\[
\left\|  \sum\nolimits_{m<n}k(f;a_{m},b_{m})w_{m}+\sum\nolimits_{N(\varepsilon
)<m}k(f;a_{m},b_{m})w_{m}\right\|  _{W}\leq2\varepsilon_{0}.
\]
So, up to $2\varepsilon_{0}$, the norm $\left\|  f\right\|  _{K}$\ of a
function $f$ may be computed by using only the vectors $\left(  w_{i}\right)
_{i=n}^{N}$ provided its support is sufficiently small. Let $\ \left(
f_{n}\right)  \subset K$ be a disjoint sequence of functions. Clearly,
$\lim_{n\rightarrow\infty}\operatorname*{mes}(\operatorname*{supp}(f_{n}))=0$.

Chose a subsequence $(f_{n_{k}})\subset\left(  f_{n}\right)  $ by an induction.

Let $n_{1}=1.$ Assume that $f_{n_{1}}$, $f_{n_{2}}$, ..., $f_{n_{m-1}}$ are
already chosen. Let $n_{k}$ be the least number $m$ such that
$s(\operatorname*{mes}(\operatorname*{supp}(f_{m}))\leq a_{j}^{-1}b_{j}$,
where $j=j(m,2^{-m}\varepsilon)$. Clearly, the sequence $(f_{n_{k}})$ spans a
subspace of $K$, which is isomorphic to a subspace of $W$ that is spanned by a
some block sequence with respect to $\left\{  w_{i}:i<\omega\right\}  $.
\end{proof}

\begin{corollary}
Let $E,F,W,\{a_{i};b_{i}:i<\omega\}$ be as in definition 1. If the space $W$
is reflexive then $K(E,F;W)$ is reflexive too. If $W$ is weakly sequentially
complete then $K(E,F;W)$ also has this property.
\end{corollary}

\begin{proof}
If $K$ is not reflexive then it contains a subspace $Y$, which is isomorphic
either to $l_{1}$ or to $c_{0}$. If $K$ is not weakly sequentially complete
then it contains a subspace $Y$, which is isomorphic to $c_{0}$. In both cases
$K$ contains a disjoint sequence of functions that spanned a subspace
isomorphic either to $l_{1}$ or to $c_{0}$. However this contradicts with
properties of $W$.
\end{proof}

\section{References}

1. Semyonov, E.M. \textit{Inclusion theorems for Banach spaces of measurable
functions}, Dokl. AN\ SSSR, \textbf{156} (1964) 1294-1295 (in Russian)

2. Plichko, A.N. and Tokarev, E.V. \textit{Bases of symmetric function
spaces}, Math. Notes, \textbf{42}, no. 1/2 (1987) 630-634 (transl. from Russian)

3. Rodin, V.A. and Semyonov, E.M. \textit{Rademacher series in symmetric
spaces,} Analysis Math., \textbf{1:4} (1975) 207-222

4. Tokarev, E.V. \textit{Subspaces of symmetric spaces of functions}, Funct.
Anal. and Appl., \textbf{13} (1979) 152-153 (transl. from Russian)

5. Tokarev, E.V. \textit{Complemented subspaces of symmetric function spaces,}
Math. Notes, \textbf{32} (1982) 926-928 (transl. from Russian)

6. Peetre, J. \textit{A new approach to interpolation spaces}, Studia Math.,
\textbf{34} (1970) 23-43

7. Beauzamy, B. \textit{Espaces interpolationes re\'{e}l; topologie et
g\'{e}om\'{e}trie,} Lect. Notes in Math. \textbf{666} (1978) 1-90

8. Davis, W.T., Figiel, T., Johnson, W.B. and Pe\l czy\'{n}ski, A.
\textit{Factoring weakly compact operators}, Journ. Funct. Anal., \textbf{17}
(1974) 317-327

9. Johnson, W.B., Maurey, B., Schechtman, G. and Tzafriri, L.
\textit{Symmetric structures in Banach spaces,} Memoirs Amer. Math. Soc.,
\textbf{217} (1979) 1-298
\end{document}